\documentclass[10pt,reqno]{amsart}

\usepackage[T1]{fontenc}
\usepackage[utf8]{inputenc}
\usepackage[english]{babel}
\usepackage{lmodern}
\usepackage{microtype}
\usepackage{amsmath,amssymb,amsthm,mathtools}
\usepackage{xcolor}
\usepackage{needspace}
\usepackage[colorlinks=true,citecolor=blue,linkcolor=blue,urlcolor=blue]{hyperref}

\newtheorem{theorem}{Theorem}
\newtheorem{corollary}[theorem]{Corollary}

\newcommand{\keyK}{\mathcal{K}}
\newcommand{\atomA}{\mathcal{A}}
\newcommand{\flagLLT}{\mathcal{G}}
\newcommand{\macH}{\mathsf{H}}
\newcommand{\macdonaldE}{\mathrm{E}}
\newcommand{\defin}[1]{%
  \relax\ifmmode\textcolor{blue}{#1}%
  \else\textcolor{blue}{\emph{#1}}\fi}

\title[Counterexamples to atom positivity]
{Counterexamples to flagged LLT and nonsymmetric Macdonald atom positivity}
\author{Per Alexandersson}
\date{}

\begin{document}

\begin{abstract}
We give a key polynomial of degree fourteen whose expansion in Demazure
$t$-atoms has a coefficient with a negative term.  This computation
disproves the conjectured Demazure atom positivity of flagged LLT
polynomials and the corresponding stable atom positivity assertion for
modified nonsymmetric Macdonald polynomials.
\end{abstract}

\maketitle

Blasiak--Haiman--Morse--Pun--Seelinger conjectured that flagged LLT
polynomials indexed by ordinary skew diagrams and the standard compatible
permutation are Demazure atom positive
\cite[Conjecture~6.6.1]{BlasiakHaimanMorsePunSeelinger2025x}.  Their construction
would also imply stable atom positivity for modified nonsymmetric Macdonald
polynomials.

\section{The key counterexample}

Let $x=(x_1,\dotsc,x_n)$ and let $s_i$ interchange $x_i$ and $x_{i+1}$,
or the corresponding entries of a composition.
The divided difference and \defin{isobaric divided difference} operators are
\[
 \defin{\partial_i f}\coloneqq\frac{f-s_i f}{x_i-x_{i+1}},
 \qquad
 \defin{\pi_i f}\coloneqq\partial_i(x_i f)
 =\frac{x_i f-x_{i+1}s_i f}{x_i-x_{i+1}}.
\]
For a weak composition $\alpha\in\mathbb{Z}_{\geq0}^n$, write $\alpha^+$ for its
weakly decreasing rearrangement, $|\alpha|=\sum_i\alpha_i$ for its size,
$\operatorname{rev}(\alpha)=(\alpha_n,\dotsc,\alpha_1)$, and
$x^\alpha=x_1^{\alpha_1}\dotsm x_n^{\alpha_n}$.
The \defin{key polynomials}, or Demazure characters, are determined by
\[
 \keyK_\alpha(x)=x^\alpha\quad\text{if }\alpha=\alpha^+,
 \qquad
 \keyK_\alpha(x)=\pi_i\keyK_{s_i\alpha}(x)
 \quad\text{if }\alpha_i<\alpha_{i+1}.
\]
The braid relations make this recursion independent of the choices of $i$.
In particular, $\keyK_{\operatorname{rev}(\lambda)}(x)=s_\lambda(x)$ for a
partition $\lambda$ padded to length $n$, so keys generalize Schur polynomials
\cite{Lascoux1990Keys}.

Replacing $\pi_i$ in this recursion by $\pi_i-1$ gives the
\defin{Demazure atoms} $\atomA_\alpha(x)$, with Mason's composition indexing
\cite[Section~2]{Mason2009}.  We similarly define
the \defin{Demazure $t$-atoms} $\atomA_\alpha(x;t)$ using
\[
 \defin{T_i}\coloneqq(1-t)(\pi_i-1)+t s_i.
\]
The initial condition is $\atomA_\lambda(x;t)=x^\lambda$ for weakly
decreasing $\lambda$.  Thus $\atomA_\alpha(x;0)=\atomA_\alpha(x)$, while
$\atomA_\alpha(x;1)=x^\alpha$.  The $t$-atoms form a basis over $\mathbb{Q}(t)$;
we write
\begin{equation}\label{eq:transition}
 \keyK_\alpha(x)=
 \sum_{\substack{\gamma\in\mathbb{Z}_{\geq0}^n\\|\gamma|=|\alpha|}}
 c_{\alpha,\gamma}(t)\atomA_\gamma(x;t).
\end{equation}
Schur polynomials have nonnegative expansions in these $t$-atoms
\cite[Corollary~31]{Alexandersson2015gbMacdonald}, so it is natural to ask
whether the same holds for all key polynomials.  The following example shows
that it does not.

\Needspace{10\baselineskip}
\begin{theorem}\label{thm:keyCounterexample}
Set
\[
 \alpha=(0,(0,1,2)^4,2)
       =(0,0,1,2,0,1,2,0,1,2,0,1,2,2),
 \qquad \gamma=(1^{14}).
\]
Then $\alpha^+=(2^5,1^4,0^5)$ and
\begin{equation}\label{eq:coefficient}
 \begin{split}
 c_{\alpha,\gamma}(t)={}&t^5+t^6+2t^7-t^{10}+5t^{11}+10t^{12}
 +26t^{13}+42t^{14}\\
 &+58t^{15}+66t^{16}+59t^{17}+42t^{18}+22t^{19}
 +7t^{20}+t^{21}.
 \end{split}
\end{equation}
In particular, its only negative term is $-t^{10}$, and
$c_{\alpha,\gamma}(1)=341$.

Moreover, for every $n\leq13$ and every
$b\in\{0,1,2\}^n$ with $|b|=n$, the coefficient
$c_{b,(1^n)}(t)$ has no negative terms.  Thus degree fourteen is minimal
within this family.
\end{theorem}

\begin{proof}
For a weak composition $b$ of length $n$ and size $n$,
nonsymmetric Hall--Littlewood duality gives
\begin{equation}\label{eq:duality}
 c_{b,(1^n)}(t)=
 [x_1\dotsm x_n]\,\keyK_b(x)
 \prod_{1\leq i<j\leq n}\frac{1-x_i/x_j}{1-tx_i/x_j}.
\end{equation}
We construct $\keyK_b$ from the defining divided differences and expand
\[
 \frac{1-z}{1-tz}=1+\sum_{r\geq1}(t-1)t^{r-1}z^r.
\]
Taking $b=\alpha$, exact coefficient extraction in \eqref{eq:duality} gives
\eqref{eq:coefficient}.

For the minimality assertion, put $h_i=i-\sum_{j=1}^i b_j$.
If some $h_i<0$, then \eqref{eq:duality} gives $c_{b,(1^n)}(t)=0$.
If $h_k=0$ for $0<k<n$, the same coefficient extraction factors as
\[
 c_{b,(1^n)}(t)=
 c_{(b_1,\dotsc,b_k),(1^k)}(t)\,
 c_{(b_{k+1},\dotsc,b_n),(1^{n-k})}(t).
\]
Indeed, every monomial of $\keyK_b$ has first $k$ exponents summing to at
least $\sum_{j=1}^k b_j=k$.  Only equality can contribute, so no factor
$x_i/x_j$ with $i\leq k<j$ contributes a positive power; the terms of the
key with equality form the product of the two smaller keys.
Since $c_{(1),(1)}(t)=1$, we need only check compositions
$b\in\{0,1,2\}^n$ with $2\leq n\leq13$, $|b|=n$, and $h_i>0$ for
$1\leq i<n$.  There are $9{,}360$ such compositions.  Exact coefficient
extraction gives no negative term in any of them.
\end{proof}

The proper prefix heights $h_i=i-\sum_{j=1}^i\alpha_j$, for $1\leq i<14$, are
$h=(1,(2,2,1)^4)$.  Thus negativity already occurs when all prefix heights
are at most two.

\section{Flagged LLT and Macdonald consequences}

LLT polynomials are generating functions of tuples of semistandard tableaux,
weighted by an inversion statistic.  Their \defin{flagged} versions
$\flagLLT_{\boldsymbol{\nu},\sigma}[X_1,\dotsc,X_\ell;t]$ are functions of
formal alphabets $X_i$, where $\boldsymbol{\nu}$ has $\ell$ rows in total.
Number the rows consecutively through the components.  Row $j$ has
\defin{flag number} $\ell+1-\sigma^{-1}(j)$.
For ordinary skew diagrams, the \defin{standard compatible permutation}
assigns flag numbers $\ell,\dotsc,1$ in increasing order of
(end content, row number), where the end content is the column index minus
the row index of the cell immediately to the right of a row
\cite[Definition~6.1.6 and Lemma~6.1.7]{BlasiakHaimanMorsePunSeelinger2025x}.
The \defin{standard flag} is the specialization $X_i=x_i$, under which
flag number $i$ gives entry bound $i$.  We use the flagged inversion
statistic and tableau formula of
\cite[Definition~6.4.8 and Corollary~6.4.11]{BlasiakHaimanMorsePunSeelinger2025x};
Weyl symmetrization recovers ordinary LLT polynomials.

The \defin{nonsymmetric Macdonald polynomials}
$\macdonaldE_\alpha(x;q,t)$ are two-parameter polynomials indexed by weak
compositions.  We use the indexing of \cite{HaglundHaimanLoehr2008}, so that
$\macdonaldE_\alpha(x;0,0)=\keyK_\alpha(x)$.
We use the \defin{modified $r$-nonsymmetric Macdonald polynomials}
$\macH_{\eta\mid\lambda}(x,Y;q,t)$ of
\cite[Definition~7.6.5]{BlasiakHaimanMorsePunSeelinger2025x}, obtained from
right-stabilized integral forms by the nonsymmetric analogue of
$f[X]\mapsto f[X/(1-t)]$.  Here $\eta$ has length $r$, $\lambda$ is a
partition, and these polynomials are nonsymmetric in $x=(x_1,\dotsc,x_r)$
and symmetric in $Y=(y_1,y_2,\dotsc)$.  Their \defin{stable atom}
basis $\atomA_{\rho\mid\lambda}(x,Y)$ is obtained by Weyl symmetrizing
ordinary atoms in the $Y$ variables and taking the stable limit
\cite[Section~7.7]{BlasiakHaimanMorsePunSeelinger2025x}.

Write
\[
 \alpha^{\mathrm{rev}}=\operatorname{rev}(\alpha)
       =(\beta,2,1,0,0),\qquad
 \beta=(2,2,1,0,2,1,0,2,1,0).
\]
We write $[B]F$ for the coefficient of a basis element $B$ in $F$.
For $n\geq1$, let $G_n(x;t)$ be the standard-flag LLT polynomial on $n$
copies of the one-box diagram $(1)$, all at the same content:
\[
 \defin{G_n(x;t)}\coloneqq
 \flagLLT_{((1),\dotsc,(1)),\mathrm{id}}[x_1,\dotsc,x_n;t]
 =\sum_{\substack{1\leq w_j\leq n+1-j\\1\leq j\leq n}}
 t^{\operatorname{inv}(w)}x_{w_1}\dotsm x_{w_n}.
\]
Here $w=(w_1,\dotsc,w_n)$ and
$\defin{\operatorname{inv}(w)}\coloneqq
|\{(i,j):1\leq i<j\leq n,\ w_i>w_j\}|$.

\begin{corollary}\label{cor:consequences}
With $c(t)=c_{\alpha,(1^{14})}(t)$ as in \eqref{eq:coefficient}, we have
\begin{align}
 [\atomA_{\alpha^{\mathrm{rev}}}(x)]G_{14}(x;t)
   &=t^{91}c(t^{-1}),\label{eq:flaggedCounterexample}\\
 [\atomA_{\alpha^{\mathrm{rev}}\mid\varnothing}(x,Y)]
 \macH_{(1^{14})\mid\varnothing}(x,Y;q,t)
   &=c(t),\label{eq:stableCounterexample}\\
 [\atomA_{\beta\mid(2,1)}(x,Y)]
 \macH_{(1^{10})\mid(1^4)}(x,Y;q,t)
   &=c(t).\label{eq:partialCounterexample}
\end{align}
Hence the first coefficient contains $-t^{81}$, while each stable-atom
coefficient contains $-q^0t^{10}$.
\end{corollary}

\begin{proof}
Corollary~6.3.6, equation~(191), of
\cite{BlasiakHaimanMorsePunSeelinger2025x}, together with key--atom duality, gives
for $n$ singleton diagrams at the same content
\begin{equation}\label{eq:translation}
 t^{\binom n2}\flagLLT_{((1),\dotsc,(1)),\mathrm{id}}
 [x_1,\dotsc,x_n;t^{-1}]
 =\sum_{\substack{b\in\mathbb{Z}_{\geq0}^n\\|b|=n}}
 c_{b,(1^n)}(t)\atomA_{\operatorname{rev}(b)}(x).
\end{equation}
The reversal in \eqref{eq:translation} is the one in the key--atom Cauchy
identity \cite[equations~(43)--(44)]{BlasiakHaimanMorsePunSeelinger2025x}.
Taking $n=14$ proves \eqref{eq:flaggedCounterexample}.

The singleton specialization of
\cite[Theorem~7.6.6(b)]{BlasiakHaimanMorsePunSeelinger2025x} is
\[
 \macH_{(1^N)\mid\varnothing}(x,Y;q,t)
 =t^{\binom N2}G_N(x;t^{-1}),
\]
which proves \eqref{eq:stableCounterexample}.  Finally,
\cite[Proposition~7.8.3]{BlasiakHaimanMorsePunSeelinger2025x} sends
$\macH_{(1^{14})\mid\varnothing}$ to
$\macH_{(1^{10})\mid(1^4)}$ after four partial Weyl symmetrizations.  The
weakly decreasing tail $(2,1,0,0)$ becomes
$\atomA_{\beta\mid(2,1)}$, which proves
\eqref{eq:partialCounterexample}.
\end{proof}

Equation~\eqref{eq:flaggedCounterexample} disproves
\cite[Conjecture~6.6.1]{BlasiakHaimanMorsePunSeelinger2025x}, while
\eqref{eq:stableCounterexample} and \eqref{eq:partialCounterexample}
disprove the stable atom-positivity assertion in Section~7.7 there.

\subsection*{AI assistance}
The counterexamples were found by ChatGPT-6 Astra.

\bibliographystyle{amsalpha}
\bibliography{bibliography}

\end{document}